\documentclass[12pt,twoside,a4paper]{amsart}
\usepackage{amssymb}
\date{\today}

\def\End{{\rm End}}

\def\deg{\text{deg}\,}

\def\sr{\stackrel}
\def\lsigma{\text{\Large $\sigma$}}

\def\w{\wedge}

\def\dbar{\bar\partial}

\def\C{{\mathbb C}}
\def\w{{\wedge}}

\def\supp{\text{supp}\,}
\def\D{{\mathcal D}}

\def\S{{\mathcal S}}

\def\Hom{{\rm Hom\, }}
\def\codim{{\rm codim\,}}

\def\Im{{\rm Im\, }}

\def\E{{\mathcal E}}

\def\O{{\mathcal O}}

\def\U{{\mathcal U}}

\def\be{\begin{equation}}
\def\ee{\end{equation}}

\newtheorem{thm}{Theorem}[section]
\newtheorem{lma}[thm]{Lemma}
\newtheorem{cor}[thm]{Corollary}
\newtheorem{prop}[thm]{Proposition}

\theoremstyle{definition}

\theoremstyle{remark}

\newtheorem{preremark}{Remark}
\newtheorem{preex}{Example}

\newenvironment{remark}{\begin{preremark}}{\qed\end{preremark}}
\newenvironment{ex}{\begin{preex}}{\qed\end{preex}}

\numberwithin{equation}{section}

\title[Integral representation with weights II]
{Integral representation with weights II, division and interpolation}

\begin{document}

\date{\today}

\author{Mats Andersson}

\address{Department of Mathematics\\Chalmers University of Technology and the University of 
G\"oteborg\\S-412 96 G\"OTEBORG\\SWEDEN}

\email{matsa@math.chalmers.se}

\subjclass{32 A 26, 32 A 27}


\thanks{The author was
  partially supported by the Swedish  Research Council}

\begin{abstract}
Let $f$ be a $r\times m$-matrix of holomorphic functions that is
generically surjective.
We provide explicit integral representation of holomorphic $\psi$
such that $\phi=f\psi$, provided that $\phi$ is holomorphic and annihilates
a certain residue current with support on the set where $f$ is not
surjective.
We also consider formulas for interpolation.
As applications we obtain generalizations
of various  results previously known for  the case  $r=1$.
\end{abstract}


\maketitle

\section{Introduction}

This paper is a continuation of \cite{A1} where we introduced a new way to generate
weighted representation formulas for holomorphic functions, generalizing \cite{BA}.
In this paper we focus on  division and interpolation
and we introduce new formulas for matrices of holomorphic functions.
As applications  we obtain generalizations
of various  results previously known for a row matrix.

Let $f=(f_1,\ldots,f_m)$ be a tuple of holomorphic functions defined in, say, a neighborhood
of the closure of the  unit ball $D$ in 
$\C^n$  with common zero set $Z$, and assume that $df_1\w\ldots\w df_n\neq 0$ on $Z$.
In \cite{BBB} was constructed a representation formula
\begin{equation}\label{rep0}
\phi(z)=f(z)\cdot \int_\zeta T(\zeta,z)\phi(\zeta) +\int_\zeta S(\zeta,z)\phi(\zeta),
\quad z\in D,
\end{equation}
for holomorphic functions $\phi$,
where both $T$ and $S$ are holomorphic in $z$, $T(\cdot, z)$ 
is  integrable,   and 
$S(\cdot,z)$ is a current of order zero (i.e., with measure coefficients)  with support on $Z$.
If $\phi$ vanishes on $Z$, thus \eqref{rep0} provides an explicit representation
of $\phi$ as an element of the ideal generated by $f$.
Moreover, if $\phi$ is just defined on $Z$, then 
$$
\int_\zeta  S(\zeta,z)\phi(\zeta)
$$
is a holomorphic extension, i.e., a holomorphic function in $D$ that interpolates
$\phi$ on $Z$.
The formula \eqref{rep0}  was extended to the case where  $f$ defines a  complete
intersection, i.e., $\codim Z=m$,   in \cite{P1}. In this case, as expected,
$S(\cdot,z)$ is closely related to the Coleff-Herrera current
$$
\dbar\frac{1}{f_1}\w\ldots\w\dbar\frac{1}{f_m};
$$
however also the singularities of $T(\cdot,z)$ are rather complicated, and contains
terms that are concentrated on $Z$, see also  Remark~\ref{sing} below. 
Formulas of this kind have  been used for various purposes by several authors;
notably for instance the explicit proof of the duality theorem
for a complete intersection in \cite{P1}, explicit versions of the fundamental
principle, \cite{BP}, sharp approximation by polynomials
\cite{Z}, and estimates of solutions to the Bezout equation, \cite{BY0};   for further examples
see  \cite{BGVY}  and the references given there. More recent
applications  can be found in \cite{DGSY} and \cite{BY}.
One can also  use  such formulas to obtain sharp estimates at the boundary, such as $H^p$-estimates,
of explicit solutions  to division problems, \cite{AC}.

In \cite{Maz} and \cite{A1} independently, was constructed a similar formula
where $T(\cdot,z)$ has quite simple,   principal value,  singularites but instead spread out over
the larger set $Y=\{f_1f_2\cdots f_m=0\}$.
In \cite{A2} we introduced a new formula like \eqref{rep0} for an arbitrary
$f$, where the singularity  of $T(\cdot,z)$ is a principal value at  $Z$,  and
$S(\cdot,z)$ is a smooth form   times a  Bochner-Martinelly type residue
current $R^f(\zeta)$ with support on $Z$.

The purpose of this paper is to extend this kind of formulas to the case
when $f$ is a generically surjective   $r\times m$ matrix of holomorphic functions. 
Given such a matrix it was defined in \cite{A5} an associated
(matrix-valued) current $R^f$ with support on the analytic set 
$$
Z=\{z;\ f(z)\ \text{is\ not\ surjective}\},
$$
with the property that if $\phi$ is an $r$-column of holomorphic
functions such that $R^f\phi=0$, then $f\psi=\phi$ has  holomorphic solutions
$\psi$ locally. In the generic case, i.e., $\codim Z=m+r-1$, 
the converse also holds. When $r=1$, and $\codim Z=m$
this is precisely the duality theorem for a complete intersection
(\cite{DS} and \cite{P1}). 
In particular we get an explicit proof of the following statement
from \cite{A5}, generalizing the Brian\c con-Skoda theorem, \cite{BS}
(for an explicit  proof in the case $r=1$, see \cite{BGVY} and \cite{EG}):
{\it  Suppose that $f$ is an  $r\times m$ matrix of holomorphic functions
that is generically surjective, and that $\phi$ is an $r$-column of holomorphic
functions. If $\|\phi\|^2\le C (\det  ff^*)^{\min(n,m-r+1)}$, then
$R\phi=0$ and hence $\phi=f\psi$  locally.}  
Here $\|\phi\|$ is the norm
$\|\phi\|^2=\langle\widetilde{ff^*}\phi,\phi\rangle$, where 
$\widetilde{ff^*}$  is the transpose of the co-matrix of $ff^*$.

We also obtain the  following result, which for the case $r=1$
appeared  in \cite{A3}. Let $\dbar^\alpha=\partial^\alpha/\partial\bar z^\alpha$
for multiindices $\alpha$.

\begin{thm}\label{co}
 Suppose that $f$ is an  $r\times m$-matrix of holomorphic functions
that is generically surjective. Let $\phi$ be an $r$-column of smooth 
functions such that 
\begin{equation}\label{villkor}
R^f (\dbar^\alpha \phi) =0
\end{equation}
for all $\alpha$. Then $\phi=f\psi$ has a smooth solution $\psi$.
In case $\codim Z=m-r+1$ the condition \eqref{villkor} is also necessary.
\end{thm}

We also present variants of this result for lower regularity.
The   division formulas admit sharp estimates at the boundary,
and as an example we indicate how one can obtain an explicit solution
of the matrix  $H^p$-corona problem.
Finally we present formulas for division problems for $\dbar$-closed forms.

\section{Representation of holomorphic functions}\label{plussig}

For a fixed point  $z$ in the open set $X$ in $\C^n$, we  let
 $\delta_{\zeta-z}$ denote interior multiplication with the vector field 
$$
2\pi i \sum (\zeta_j-z_j)\frac{\partial}{\partial \zeta_j},
$$
and let $\nabla_{\zeta-z}=\delta_{\zeta-z}-\dbar$. 
We begin with a slight  generalization of the main result in \cite{A1}
(lower indices denote bidegree).

\begin{prop}\label{vikt}
Assume that $z$ is a fixed point in $X$ and 
 $g=g_{0,0}+\ldots+g_{n,n}$ is a current  in $X$  with compact support
such that  $\nabla_{\zeta-z}g=0$. Moreover, assume that $g$ is smooth in a neighborhood
of $z$    and $g_{0,0}(z)=1$.
Then
\begin{equation}\label{rep}
\phi(z)=\int g\phi=\int g_{n,n}\phi
\end{equation}
for each holomorphic function $\phi$ in $X$.
\end{prop}

For the reader's convenience we supply the simple  proof.

\begin{proof}
Let $u=u_{1,0}+\ldots +u_{n,n-1}$ be a current that is smooth outside the point $z$ and 
such that $\nabla_{\zeta-z} u=1-[z]$, where $[z]$ denotes the $(n,n)$-current point evaluation
at $z$. For instance one can take 
$$
u=\frac{b}{\nabla_{\zeta-z}b}=b +b\w \dbar b+\ldots + b\w(\dbar b)^{n-1},
$$
where $b=\partial|\zeta-z|^2/|\zeta-z|^2 2\pi i$, see \cite{A1}.
Then $u\w g$ is a well-defined current with compact support, and 
$$
\nabla_{\zeta-z}(u\w g)=g-[z]\w g= g-[z]
$$
since $g_{0,0}(z)=1$.  Therefore,  
$\dbar (u\w g)_{n,n-1}\phi=\phi(z)[z]-g_{n,n}\phi$,
 which implies \eqref{rep} by Stokes' theorem.
\end{proof}

A form $g=g_{0,0}+g_{1,1}+\cdots +g_{n,n}$ which is smooth in a neighborhood  of our fixed point $z$,
and such that $g_{0,0}(z)=1$, will be called a weight (with respect to $z$). Notice that if 
$g^1$ and $g^2$ are weights with disjoint singular supports, then
again the  product  $g=g^1\w g^2$ is a weight.
Moreover, if  $\phi$ takes values in the vector bundle $Q\to X$, and
$g$ takes values in $\Hom(Q,Q)$, $\nabla_{\zeta-z}=0$ and $g_{0,0}(z)=I_Q$, $I_Q$ denoting the identity
morphism $Q\to Q$,   then \eqref{rep} still
holds; this follows by  the same proof as for  the scalar-valued case.

\begin{ex}\label{exett}
Assume that $D$ is a smoothly bounded domain in $X$ that admits a smooth family of holomorphic support
functions, i.e.,  $\Gamma(\zeta,z)\in C^{\infty}(\partial D\times U)$, where $U\supset\overline D$,
depending holomorphically on $z$, such that $\Gamma(\zeta,z)$ is non-vanishing for
$z\in \overline D\setminus\{\zeta\}$ and $\Gamma(z,z)=0$ for $z\in \partial D$. 
Then $D$ is necessarily pseudoconvex and  we may assume that 
$\Gamma(\zeta,z)=\delta_{\zeta-z}\gamma(\zeta,z)$, where $\gamma$ is a smooth $(1,0)$-form
that is holomorphic for $z\in D$.  If 
$s=\gamma(\zeta,z)/ \delta_{\zeta-z}\gamma(\zeta,z)$,
then for each $z\in D$,
$$
g=\chi_D-\dbar\chi_D\w\frac{s}{\nabla_\zeta s}=
1-\dbar\chi_D\w [s+s\w\dbar s+ s\w(\dbar s)^2+\cdots +s\w(\dbar s)^{n-1}]
$$
is a  weight (with respect to $z$)  with support on $\overline D$, smooth outside $\partial D$,
and depending holomorphically on $z$.
Let  $g'$ be  any weight (with respect to  $z$) that is smooth in a neighborhood of $\partial D$.  
Since $\dbar \chi_D=-[\partial D]_{0,1}$ we get from \eqref{rep}
the formula
\begin{equation}\label{r2}
\phi(z)= \int_D g' \phi +\int_{\partial D} g'\phi\w [s+s\w\dbar s+ s\w(\dbar s)^2+\cdots +s\w(\dbar s)^{n-1}].
\end{equation}
The existence of such families of holomorphic support functions for strictly pseudoconvex domains
is due to Henkin and Ramirez, see, e.g., \cite{HL}. In \cite{DF1} and \cite{DF2}
are    constructed  families  of holomorphic support functions,  admitting sharp estimates, 
for (linearly) convex domains of finite type.
If $D$ is the unit ball in $\C^n$  we can take
$s(\zeta,z)=\partial|\zeta|^2/2\pi i (1-\bar\zeta\cdot z)$;
we then  get \eqref{r2} with
$$
s\w (\dbar s)^{k-1}=\frac{1}{(2\pi i)^k}\frac{\partial|\zeta|^2\w(\dbar\partial|\zeta|^2)^{k-1}}
{(1-\bar\zeta\cdot z)^k}.
$$
\end{ex}

For our purposes it is  convenient to have a formula like this for a weight  $g'$ that is
not necessarily smooth on $\partial D$.

\begin{ex}\label{extva}
Assume that $X$ is pseudoconvex  and  let $K\subset X$ be a
holomorphically convex compact subset. Moreover let $\chi$ be a  cutoff function  that is identically
$1$ in a neighborhood of  $K$. It is easy to  find a $(1,0)$-form $s(\zeta,z)$ on the support of $\dbar\chi$,
depending holomorphically on $z$ in a neighborhood of $K$, such that $\delta_{\zeta-z}s=1$.
Then for each $z\in K$, 
$$
g=\chi-\dbar\chi\w\frac{s}{\nabla_{\zeta-z} s}=
\chi-\dbar\chi\w [s+s\w\dbar s+ s\w(\dbar s)^2+\cdots +s\w(\dbar s)^{n-1}]
$$
is a compactly supported weight that depends holomorphically on $z$.
If $K$ is the closure of the unit ball $D$ we can take
$$
s(\zeta,z)=\frac{\partial|\zeta|^2}{2\pi i(|\zeta|^2-\bar\zeta\cdot z)}.
$$
If $g'$ is any weight in $X$  (with respect to $z\in K$), then we get the representation formula
\begin{multline}\label{r1}
\phi(z)= \int \chi g'\phi -\\
\int \dbar\chi\w [s+s\w\dbar s+ s\w(\dbar s)^2+\cdots +s\w(\dbar s)^{n-1}]
\w g'\phi, \quad  z\in K.
\end{multline}
\end{ex}

\section{Division formulas in the case $r=1$}\label{case1}

To begin with, let $f$ be a row matrix  of holomorphic functions in a pseudoconvex
domain $X\subset\C^n$. In \cite{A2}
were introduced formulas for division and interpolation, and more generally,
homotopy formulas for the Koszul complex induced by $f$. In this section we
derive these formulas in a new way that will model the construction when
$r>1$. It is convenient to introduce a trivial rank $m$ bundle $E$ over $X$
and think of $f$ as a section of the dual bundle $E^*$.
If we let $\delta_f\colon\Lambda^{k+1}E\to\Lambda^k E$ denote 
 interior multiplication with $f$, we have the Koszul complex
\begin{equation}
0\to \Lambda^m E\stackrel{\delta_f}{\to} \cdots \stackrel{\delta_f}{\to}
\Lambda^2 E_2\stackrel{\delta_f}{\to}  E\stackrel{\delta_f}{\to} \C\to 0,
\end{equation}
where $\C$ is the trivial line bundle. We  consider currents with values in $\Lambda E$ as
sections of the bundle $\Lambda(E\oplus T^*(X))$, so that, e.g., 
differentials and sections of $E$ anti-commute, and 
$\delta_f$ and $\dbar$ anti-commute; for more details, see \cite{A2}.
Assume that we have  $(0,k-1)$-currents $U_k$, smooth outside some analytic variety,  and
$(0,k)$-currents $R_k$ with values in $\Lambda^k E$,
 $R_k$ having support  on $Z=\{f=0\}$,  such that 
\begin{equation}\label{fot}
(\delta_f-\dbar) U=1-R,
\end{equation}
where $U=U_1+\cdots +U_m$ and $R=R_1+\cdots +R_m$. Specific choices will be discussed below.
Assume that $\phi$ is a holomorphic section of $\Lambda^\ell E$   such that
$\delta_f\phi=0$ and 
$R\w \phi=0$. Then  it follows from \eqref{fot} that $(\delta_f-\dbar)(R\w \phi)=0$, and
by solving a sequence of $\dbar$-equations one finds (locally)  
a holomorphic section $\psi$ of $\Lambda^{\ell+1}E$  such that  $\delta_f\psi=\phi$, see \cite{A2}. 
We will now provide an explicit formula for such a solution $\psi$.

\smallskip
Let $e_j$ be a global frame for $E$ with dual frame  $e_j^*$ for $E^*$ so that 
$f=\sum _1^m f_j e_j^*$.  One can find holomorphic  $(1,0)$-forms 
$h_j$ such that $\delta_{\zeta-z}h_j=f_j(\zeta)-f_j(z)$,
so called Hefer forms. Now $h=\sum h_j\w e_j^*$ 
induces a mapping $\delta_h$, taking a $(p,q)$-current-valued section of
$\Lambda^{k+1}E$ to a $(p+1,q)$-current-valued section of $\Lambda^k E$.
Since $h$ has total degree $2$, $\delta_h$ commutes with $\delta_f$ and $\delta_{\zeta-z}$.
If  $(\delta_h)_k=\delta_h^k/k!$ and $\delta_{f(z)}$ is 
interior multiplication with the section $f(z)=\sum f_j(z) e_j^*$ of $E^*$, then
for a $(0,q)$-current $\xi$ with values in $\Lambda^kE$ we have
\begin{equation}\label{fot2}
\delta_{\zeta-z}(\delta_h)_k\xi=(\delta_h)_{k-1}(\delta_f-\delta_{f(z)})\xi
\end{equation}
for all integers $k$, if $(\delta_h)_k$ is interpreted as $1$ for $k=0$ and  zero for $k<0$.
Assume that  $\phi$ is holomorphic and takes
 values in $\Lambda^\ell E$.  Using \eqref{fot} and \eqref{fot2},
a straight forward computation shows that 
$$
g'=\delta_{f(z)}\sum_k(\delta_h)_{k-1}( U_k\w\phi) +
\sum_k (\delta_h)_{k-1}(U_k\w\delta_f\phi)
+\sum_k (\delta_h)_k (R_k\wedge \phi)
$$
is $\nabla_{\zeta-z}$-closed for each fixed $z$.
Moreover, by \eqref{fot}, $\delta_fU_1\phi=\phi$, so 
$$
g'_{0,0}(z)=\delta_{f(z)}(U_1\w\phi)|_{\zeta=z}=\phi(z).
$$
For each $z$ outside $Z$ and the set where $U$ is not smooth, by (an immediate consequence of) 
Proposition~\ref{vikt} we get the representation
\begin{equation}\label{ops}
\phi(z)=\delta_{f(z)} T\phi(z)+T(\delta_f\phi)(z) +S\phi(z),
\end{equation}
where
$$
T\phi(z)=\int
\sum_k(\delta_h)_{k-1}( U_k\w\phi)\w g, 
$$
and
\begin{equation}\label{Sformel}
S\phi(z)=\int \sum_{k} (\delta_h)_k (R_k\w\phi)\w g, 
\end{equation}
if $g$ is a smooth weight with compact support.
Since each term in \eqref{ops} is holomorphic in $z$, the equality must hold
everywhere.

\begin{ex}\label{standex}
With the following choice of currents $U$ and $R$, the representation
\eqref{ops} is precisely the formula in Theorem~9.3 in \cite{A2} expressed in a new way.
Assume that $E$ is equipped with a Hermitian metric  and let
$\sigma$ be the section of $E$ over $X\setminus Z$ with pointwise minimal
norm such that $f\sigma=1$. If 
$E$ has the trivial metric  with respect to the  global holomorphic frame $e_j$,
then 
$
\sigma=\sum_j \bar f_j e_j/|f|^2.
$
Let 
$$
u^f=\frac{\sigma}{\nabla_f\sigma}=\sigma+\sigma\w\dbar\sigma+\sigma\w(\dbar\sigma)^2+\cdots.
$$
Then the principal value current
$
U^f=\lim_{\epsilon\to 0}\chi_{|f|>\epsilon} u^f,
$
exists and is a current extension of $u^f$ across $Z$. Moreover
$
R^f=\lim_{\epsilon\to 0} \dbar \chi_{|f|>\epsilon}\w u^f,
$
exists and \eqref{fot} holds,
see \cite{A2}. Alternatively, $U^f$ and $R^f$ can be defined as
the value at $\lambda=0$ of the analytic continuation of
$|f|^{2\lambda}u^f$ and $\dbar|f|^{2\lambda}\w u^f$, respectively.
Moreover,  these  currents can  be obtained as limits of smooth forms.
Notice that the section  $s=|f|^2\sigma$ is smooth. Now, see \cite{HS1} and \cite{HS2},
$$
U^f=\lim_{\epsilon\to 0}\sum_\ell\frac{s\w(\dbar s)^{\ell-1}}
{(|f|^2+\epsilon)^{\ell}}, 
$$
and
$$
R^f=\lim_{\epsilon\to 0}\frac{\epsilon}{|f|^2+\epsilon}\sum_\ell
\frac{(\dbar s)^{\ell}}
{(|f|^2+\epsilon)^{\ell}}.
$$
For degree reasons $R^f_k=0$ if $k>\min(m,n)$ and $U^f_k=0$ if $k>\min(m,n+1)$. 
It turns out also that  $R_k^f=0$ if $k<p=\codim Z$, 
and therefore  $R^f=R^f_p+\cdots +R^f_{\min(m,n)}$.
Thus the sum in \eqref{Sformel}  only runs from
$k=p$ to $k=\min(n,m-\ell)$. 
Therefore,  if $\phi$ is a holomorphic section of $\Lambda^\ell E$ such that
$\delta_f\phi=0$ and,  in addition, $R^f\w \phi=0$, then 
$$
\psi(z)=\int_\zeta \sum_{k=1}^{\min(n+1,m-\ell)} (\delta_h)_{k-1}(U^f_k\w\phi)\w g
$$
is a holomorphic solution to $\delta_f\psi =\phi$.
\end{ex}

\begin{ex}
Assume that $f$ defines a complete intersection, and let 
$$
U_1= \frac{e_1}{f_1}, \quad  U_{k+1}=\frac{e_{k+1}}{f_{k+1}}\w\dbar U_k, \quad
R=\dbar\frac{1}{f_1}\w\ldots\w\dbar\frac{1}{f_m}\w e_m\w\ldots\w e_1,
$$
where the currents are of Coleff-Herrera type.
Then, using the calculus from  \cite{P}, cf., also   Section~5 in \cite{A1}, 
one can check that \eqref{fot} holds. Thus \eqref{ops} gives a
division formula in \eqref{ops}, which is singular over the set $\{f_1\cdots f_m=0\}$.
This formula is  similar to, but even simpler than,   the formula  in Section~5 of \cite{A1}.
\end{ex}

\begin{remark}[Berndtsson's division formula]\label{sing}
In our notation, Berndtsson's classical division formula, \cite{BBB}, 
  can be described in the following way. For $\epsilon>0$, let
$\sigma^{\epsilon}=s/(|f|^2+\epsilon)$
(here $s$ is as in Example~\ref{standex}) and let $h$ be a Hefer form as above. Then
$$
g'= 1-\nabla_\eta h\cdot\sigma^\epsilon
=\frac{\epsilon}{|f|^2+\epsilon}+f(z)\cdot\sigma^\epsilon+h\cdot\dbar\sigma^\epsilon
$$
is a weight, so  by Proposition~\ref{vikt} we have the  representation
$$
\phi(z)=\int \Big(\frac{\epsilon}{|f|^2+\epsilon}+f(z)\cdot\sigma^\epsilon+h\cdot\dbar\sigma^\epsilon
\Big)^{\min(n+1,m)}\w g \phi.
$$
if $\phi$ is a holomorphic  function (i.e.,  $\ell=0$).
Possibly besides the choice of form $g$, this is precisely the formula introduced  in \cite{BBB}.
One can prove that it  converges to a decomposition like \eqref{rep0}, 
for a quite arbitrary tuple $f$,  when $\epsilon\to 0$; 
the non-complete intersection case
is studied in \cite{BGVY} and \cite{BY} (but using analytic continuation) and
in \cite{EG}.  Making the most natural decomposition,
letting all terms without  a factor $f(z)$ together constitute $S\phi$,
the resulting current in $T\phi$  is not of simple principal value type but will involve 
terms concentrated on $Z$. However, in the case when $f$ defines
a complete intersection it seems that these bad terms disappear.
Moreover, in the the general case, 
and under the hypothesis that $|\phi|\le C|f|^{\min(m,n)}$  that is considered 
in \cite{BGVY} and \cite{EG} to get an   explicit proof of the Brian\c con-Skoda theorem
(in \cite{EG} even a more general form of this theorem is considered),
these terms vanish. Therefore it seems that with no essential loss,
one could incorperate these bad terms in $S\phi$ from the beginning.
\end{remark}

\section{Generalized Hefer forms in the unit ball}\label{heferball}

We  shall now make an explicit computation of the formula \eqref{ops}
in the case when $f$ is the tuple
$$
f(\zeta)=\sum_1^n(\zeta_j-w_j) e_j^*=(\zeta-w)\cdot e^*
$$
in the unit ball $D$ for a fixed $w\in D$; a similar computation works in any domain
that admits a smooth family of holomorphic support functions, cf., Example~\ref{exett} above.
In order to get  holomorphic dependence of $w$ we define 
$u=\sigma+\sigma\w(\dbar\sigma)+ \sigma\w(\dbar\sigma)^2+\ldots$
with
$$
\sigma=\frac{\sum\bar\zeta_j e_j}{|\zeta|^2-\bar\zeta\cdot w}=\frac{\bar\zeta\cdot e}{|\zeta|^2-\bar\zeta\cdot w}
$$
outside the singularity $w$, and since $u$ so defined is integrable we can let $U$ be the 
trivial extension across $w$. 
Since $f(\zeta)-f(z)=(\zeta-z)\cdot e^*$ we can take
$$
h=\frac{1}{2\pi i}\sum_1^n d\zeta_j\w e_j^*
$$
as our Hefer form.  Then,   see, e.g., \cite{A1},   Proposition~2.2 (just replace $e_j$ by $d\zeta_j$),
we have that $(\delta_f-\dbar)U=1-R$, where 
\begin{equation}\label{dumle}
(\delta_h)_k R_k=0, \quad k<n,\quad (\delta_h)_n R_n=[w].
\end{equation}
Since $U$ and $R$ have no singularities  at the boundary we can use the weight 
$g$ in \eqref{r2}.

We first consider \eqref{ops} when  $\phi$ is a function. 
Since   $\sum_1^n (\delta_h)_{k-1} U_k$ has no component of bidegree $(n,n)$,
formula \eqref{ops} becomes
$$
\phi(z)-\phi(w)=\delta_{f(z)}
\int_{\partial D}\sum_{k=1}^n (\delta_h)_{k-1}(\sigma\w (\dbar \sigma)^{k-1}\phi)\w s\w(\dbar s)^{n-k}.
$$
Noting that $\delta_h \sigma\w s=0$, 
more explicitly we have
$$
\phi(z)-\phi(w)=\delta_{f(z)}
\int_{\partial D} \sum_{k=1}^n\frac{\bar\zeta\cdot e\w\partial|\zeta|^2
\w (\frac{i}{2\pi}\partial\dbar|\zeta|^2)^{n-1}}
{(1-\bar\zeta\cdot w)^{k}(1-\bar\zeta\cdot z)^{n-k+1}}.
$$
Since $\phi$ clearly depends holomorphically on $w$,  this is
an explicit Hefer decomposition of the function $\phi$ in $D$.
In fact this is precisely what we  get if we express  $\phi(z)-\phi(w)$ by means of
the Szeg\H{o} integral as, e.g., in   \cite{AC}.
The new  interesting case here is when $\phi$ takes values in $\Lambda^\ell E$, $\ell>0$.
In view of \eqref{dumle}, then  $R\w  \phi=0$,  so we get instead 
$$
\phi=\delta_f T\phi +T(\delta_f \phi),
$$
where
$$
T\phi(z) =\int_{\partial D} 
\sum_{k=1}^{n-\ell}
\frac{\bar\zeta\cdot e\w(\delta_h)_{k-1}\big[(d\bar\zeta\cdot e)^k\w\phi\big]\w 
\partial|\zeta|^2 \w(\frac{i}{2\pi}\partial\dbar|\zeta|^2)^{n-k}}
{(1-\bar\zeta\cdot w)^{k}(1-\bar\zeta\cdot z)^{n-k+1}}.
$$
In particular, if $\delta_f\phi=0$ we have that
$
\delta_fT\phi = \phi.
$
It is clear from this formula that $T\phi$ depends holomorphically on $w$.
Moreover, it is proved in \cite{AC} that mappings like 
$\phi\mapsto T\phi$ admit certain sharp estimates at the boundary.

\section{Division formulas in the case $r>1$}

In order to generalize the formulas in Section~\ref{case1} to the case with matrices
$f$ we first consider a quite abstract setting. Assume that we have a finite 
 complex of Hermitian holomorphic vector bundles over $X$ 
\begin{equation}\label{kollo}
0\to E_N\sr{f_N}{\to}\ldots \sr{f_3}{\to} E_2 \sr{f_2}{\to} E_1\sr{f_1}{\to} E_0\to 0. 
\end{equation}
We will consider currents with values in $E=\oplus E_k$ and in
$\Hom(E_0,E)$, i.e.,  sections of the bundles
$\D'_{\bullet}(X,E)=\D'_\bullet(X)\otimes_{\E(X)}\E(X,E)$ and 
$\D'_{\bullet}(X,\End E)=\D'_\bullet(X)\otimes_{\E(X)}\E(X,\Hom(E_0,E))$.
Clearly, $f=\sum f_k$ and $\dbar$ act on 
these spaces and we will arrange so that $f\dbar=-\dbar f$.
To obtain this, it is natural 
to  consider $E$ as a superbundle, $E=E^+\oplus E^-$,
with $E^+=\oplus E_{2k}$ and $E^-=\oplus E_{2k+1}$,
so   that sections of $E^+$ have even degree and sections of $E^-$ have odd degree.
The space $\D'_\bullet(X,E)$ 
has a natural structure as a left $\E_{\bullet}(X)$-module, and it gets  a natural  $Z_2$-grading by
combining that gradings of $\D_\bullet(X)$ and $\E(X,E)$. We make
$\D'_\bullet(X,E)$ into a right  $\E_{\bullet}(X)$-module  by letting
$\xi \phi= (-1)^{\deg \xi \deg \phi}\phi \xi$ for sections $\xi$ of $\E(X,E)$
and smooth forms $\phi$. 
The superstructure on $E$ induces a superstructure  $\End E=\End(E)^+\oplus\End(E)^-$
so  that a mapping is odd if, like $f$,  it maps  $E^+\to E^-$ and $E^-\to E^+$. 
In the same way we get a $Z_2$-grading of $\D'_\bullet(X,\End E)$. 
For instance,  $\dbar$ extends to  an odd  mapping on $\D'_\bullet(X,E)$,
as well as on $\D'_\bullet(X,\End E)$.
Since $f$ is holomorphic and of odd degree, we have that $\dbar\circ f=-f\dbar$.

\smallskip

Let us now assume that  we have $(0,k-1)$-currents $U_k$  and
$(0,k)$-currents $R_k$ with values in $\Hom(E_0,E_k)$ such that
\begin{equation}\label{puma}
f_1 U_1=I_{E_0},\quad  f_{k+1} U_{k+1}-\dbar U_k=R_k.
\end{equation}
Moreover, we assume that $U=\sum U_k$ is smooth outside some analytic variety $Z$ and
that $R=\sum R_k$ has its support on $Z$. A possible choice of such
currents will be discussed in the next section. We will also use the short-hand notation
$(f-\dbar)  U=I_{E_0}-R$  for \eqref{puma}.    Notice that 
$f-\dbar$ is (minus) the $(0,1)$-part of the super connection $D-f$ introduced
by Quillen, \cite{Qu},  where $D$ is the Chern connection  on $E$.

\begin{prop}\label{locke}
Assume that $\phi$ is a holomorphic section of $E_0$ such that
$R\phi=0$. Then, locally,  $f_1\psi=\phi$ has  holomorphic solutions $\psi$.
\end{prop}

\begin{proof} 
In fact,  by assumption $(f-\dbar)(U\phi)=\phi$, and hence
by successively solving the equations 
$\dbar w_k=U_k\phi +f_{k+1} w_{k+1},$ we get the holomorphic solution $\psi=U_1\phi +f_2 w_2$.
\end{proof}

Our aim now is to construct an explicit formula that provides the desired solution $\psi$.
If we restrict our  attention to some open domain $X$ where we can choose
global frames for all the bundles, 
then,  for  each fixed $z\in\Omega$, $f_k(z)\colon E_k\to E_{k-1}$ is a
well-defined morphism, which coincides with $f_k$ on the fiber over $z$.
For fixed $z\in X$, as before, let $\delta_{\zeta-z}$ and $\nabla_{\zeta-z}$
be as  in Section~\ref{plussig}.

\begin{lma}\label{hlemma}
 (i)  Assume that $X$ is pseudoconvex.
For any holomorphic function $\phi$ we can find a holomorphic $(1,0)$-form $h$,
depending holomorphically on $z$,  such that
$\delta_{\zeta-z} h=\phi(\zeta)-\phi(z)$.
\smallskip

\noindent (ii)  If $\xi$ is a holomorphic $(k,0)$-form, $k\ge 1$, depending holomorphically on the 
parameter $z$, such that $\delta_{\zeta-z} \xi=0$, then we can find a holomorphic  
$(k+1,0)$-form $\xi'$ depending holomorphically on $z$ such that 
$\delta_{\zeta-z}\xi'=\xi$. 
\end{lma}

These facts are well-known and follow from  Cartan's theorem.
For an explicit construction in the unit ball, see  Section~\ref{heferball} above.

\begin{prop}[Existence of Hefer forms]\label{hprop}
There are $(k-\ell,0)$-form-valued holomorphic morphisms $H^\ell_k\colon E_k\to E_\ell$,
depending holomorphically on $z$, such that 
$H^\ell_k=0$ for $ k<\ell$, $ H^\ell_\ell=I_{E_\ell}$, and in general,
\begin{equation}\label{Hdef}
\delta_{\zeta-z} H^\ell_{k}=
H^\ell_{k-1}f_k(\zeta) -f_{\ell+1}(z) H^{\ell+1}_{k}.
\end{equation}
\end{prop}

\begin{proof}
In fact, for $k=\ell+1$ the right hand side of \eqref{Hdef} is just
$f_{\ell+1}(\zeta)-f_{\ell+1}(z)$ so the existence of $H^\ell_{\ell+1}$ is
ensured by the first part of the lemma applied to the entries in the  matrix
repesentation of  $f_{\ell+1}$. If $k>\ell+1$, then the right hand
side of \eqref{Hdef} is $\delta_{\zeta-z}$-closed, so in view of the second part of the
lemma, the proposition follows by induction over $\ell$ down-wards,
starting with $\ell=N-1$, and over $k$ up-wards.
\end{proof}

Assuming $U$ is smooth outside $Z$ and $R$ supported on $Z$, for fixed $z\notin Z$,  we can define
the current
$$
g'=f_1(z) \sum_1^\mu H^1_k U_k+ \sum_1^\mu H^0_k R_k,
$$
and it is easily checked that 
\begin{equation}\label{omsk}
\nabla_{\zeta-z} g'= 0, \quad\quad  g'_{0,0}(z)=I_{E_0}.
\end{equation}
In fact, noticing that \eqref{Hdef} holds for all $k$ and $\ell$,
we can write
$g'=f(z) H^1 U +H^0 R$  and use \eqref{Hdef} to get 
\begin{multline*}
\nabla_{\zeta-z} g'=(\delta_{\zeta-z}-\dbar)g'= \\
-f(z)\Big((H^1\delta-f(z)H^2)U-H\dbar U\Big)+
(H^0f(\zeta)-f(z) H^1)R- H^0\dbar R= \\
-f(z) H^1  I_{E_0} =0,
\end{multline*}
where we have used that $f\circ f=0$ (this holds since \eqref{kollo} is a complex) and
that $f=0$ on $E_0$.

\begin{prop}\label{skrutt}
Let $\phi$ be a holomorphic section of $E_0$ in $X$,
let $g$ be a scalar weight with compact support
in $X$ for each $z$ in  $\U \subset\subset X$,  and assume that $g$ 
depends holomorphically on $z$. Then we have the holomorphic decomposition
\begin{equation}\label{formel}
\phi(z)=
f_1(z)\int_\zeta H^1U \phi\w g + \int_\zeta H^0 R \phi \w g,\quad z\in \U.
\end{equation}
\end{prop}

\begin{proof}
In view of \eqref{omsk} the case when $z\notin Z$ follows from Proposition~\ref{vikt}. Since each term 
in  \eqref{formel} is  holomorphic for $z$ in $\U$, the equality must hold
also  across $Z$.
\end{proof}

In particular we see directly that 
\begin{equation}\label{soll}
T\phi=\int H^1 U \phi\w g
\end{equation}
is a holomorphic solution to $f_1\psi=\phi$ in   $\U$ if
$R\phi=0$. We have thus obtained  an explicit representation
of the solution in Proposition~\ref{locke}.

\smallskip

One can also use \eqref{formel} for interpolation.

\begin{prop}
Assume that $\phi$ is a holomorphic section of $E_0$ in a neighborhood  of $Z$ in $X$.  Then 
$$
S\phi=\int_\zeta H^0R\phi\ g
$$
is a holomorphic section $E_0$  in $\U$ such that
$\phi-S\phi$ belongs to the image of $f_1$ locally at $Z\cap \U$.
\end{prop}

\begin{proof}
Recall that  $S\phi$ only depends on $R\phi$ and thus only depends on
the values of $\phi$ on $Z$ up to the order of the current $R$.
In virtue of Cartan's theorem there is some holomorphic section $\Phi$
in $X$ that coincides with $\phi$ up to the prescribed order on $Z$.
From \eqref{formel} it then follows that
$S\phi=S\Phi=\Phi-f_1 T\Phi$ from which the proposition follows.
\end{proof}

It is possible to give a direct argument of this proposition, with no
reference to Cartan's theorem, cf., Remark~3 in \cite{A2}.
In fact, suppose that $\phi$ is holomorphic in the open set
$U\supset Z$, and  take a cutoff function $\chi$ with support in
$U$ and equal to $1$ in a small neighborhood of $Z$.
For a fixed $z^0$ on $Z$ we can 
find a $(1,0)$-form $s(\zeta)$ such that $\delta_{\zeta-z^0}s(\zeta)\neq 0$ for $\zeta$ 
on $\supp\dbar\chi$. 
By continuity this will hold also for $z$ in a small neighborhood $V$ of $z^0$.
Therefore, $g'=\chi-\dbar\chi\w s/\nabla_{\zeta-z}s$
is a weight for each $z\in V$ and hence
\begin{equation}\label{apa}
\phi(z)=f_1(z)\int H^1 U\phi\w g\w g'+\int H^0 R\phi \w g\w g', \quad z\in V.
\end{equation}
However, since $g'\equiv 1$ in a neighborhood of $Z$, the last term
coincides with $S\phi(z)$ for $z\in V$, and hence \eqref{apa} shows that  
$\phi-S\phi$ is in the image of $f_1$ there.


\section{Generically surjective morphisms}

Our main application in this paper of the abstract procedure developed in the
preceding section is when $E$ and $Q$ are given
(trivial) Hermitian holomorphic vector bundles and $f\colon E\to Q$ is a pointwise surjective,
or possibly generically surjective, holomorphic morphism;  we look for an explicit
formula for a holomorphic solution to $f\psi=\phi$, where $\phi$ is a section of $Q$.
Let us assume that $E$ and $Q$ are bundles over a neighborhood
of the closed unit ball and let  $Z$ be the analytic variety where $f$ is not surjective. 
Following \cite{A5} we  take  $E_0=Q$, $E_1=E$, and 
$$
E_k=\Lambda^{k+r-1}E\otimes S^{k-2} Q^*\otimes \det Q^*,  \quad k\ge 2.
$$
Let $e_j$ and $\epsilon_k$ be global  holomorphic frames for $E$ and $Q$
respectively, and let $e_j^*$ and $\epsilon_k^*$ be their dual frames.
Then   $f=\sum f^k\otimes\epsilon_k$, where $f^k$ are sections of $E^*$
and $\det f=f^1\w\ldots\w f^r\otimes\epsilon_r\w\ldots\w\epsilon_1$.
We obtain a complex \eqref{kollo} by taking 
$f_1=f$, $f_2=\det f$ and $f_k$ as interior multiplication  $\delta_f$ 
with $f$ for $k>2$, see \cite{A5} for details, which is known as the   Eagon-Northcott complex. 
In the case $r=1$ this is just the Koszul complex.
If $r$ is odd, the natural grading of $\Lambda(T^*_{0,1}(X)\oplus E)$ gives rise to the
desired $Z_2$-grading of $\oplus E_j$; if $r$ is even, then $\dbar$ and $\det f$ do not
anti-commute and therefore one has to compensate with a factor $(-1)^{r+1}$ at some places,
see \cite{A5}; in what follows it is therefore tacitly understood that
$r$ is odd, and we leave it to the interested reader to find out where to put  necessary minus signs
in the case when $r$ is even.

Outside $Z$ let  $\sigma_k$ be the sections of $E$ with  minimal norms such that
$f^j\sigma_k=\delta_{jk}$. Then $\sigma=\sigma_1\otimes\epsilon_1^*+\ldots+\sigma_r\otimes\epsilon_r$ 
is the minimal section of $\Hom(E_0,E_1)$ such that $f\sigma=I_{E_0}$.
Moreover, the section $\det f$  is nonvanishing, and if
$\lsigma=\sigma_1\w\ldots\w\sigma_r\otimes\epsilon_r^*\w\ldots\w\epsilon_1^*$,
for $\xi$ with values in $E_1$, then $\lsigma\xi$ is the minimal inverse.
Now we  define 
\begin{equation}\label{udef1}
u_1=\sigma, \quad 
u_k=(\dbar\sigma)^{\otimes(k-2)}\otimes\lsigma\otimes\dbar\sigma, \quad k\ge 2,
\end{equation}
where $\otimes$ shall be interpreted as  $\w$ on the factors in $\Lambda(T^*_{0,1}(X)\oplus E)$
and $\otimes$ on the $Q^*$- factors, and
where the rightmost  factor $\dbar\sigma$ is supposed to act on 
$E_0$. If we 
extend $u$ across $Z$ as $U=|\det f|^{2\lambda} u|_{\lambda=0}$, cf., Example~\ref{standex}
and see \cite{A5} for details, and  let $R=\dbar|\det f|^{2\lambda}\w u|_{\lambda=0}$,
then we get currents satisfying \eqref{puma}.
Moreover, $U$ is smooth outside $Z$ and $R$ has support on  $Z$.

To construct the division formula we also need suitable Hefer forms, whose existence
are ensured by Proposition~\ref{hprop}. However, we can  be
somewhat more explicit.
To begin with let  $h(\zeta,z)$ be a $(1,0)$-form with values in
$\Hom(E_1,E_0)$ such that $\delta_{\zeta-z}h=f(\zeta)-f(z)$, and let
$(\delta_h)_\ell=(\delta_h)^\ell/\ell!$. Notice that since
$\delta_f$ is an odd mapping, $\delta_h$ is even.
For  $\ell\ge 2$ we can now take  $H^\ell_k=(\delta_h)_{k-\ell}$. In fact,
$$
\delta_{\zeta-z}(\delta_h)_{k-\ell}=(\delta_h)_{k-\ell-1}(\delta_{f(\zeta)}-\delta_{f(z)})=
(\delta_h)_{k-\ell-1}\delta_{f(\zeta)}-\delta_{f(z)}(\delta_h)_{k-\ell-1},
$$
which shows that \eqref{Hdef} is fulfilled for $\ell\ge 2$.

\begin{thm}
Let $f\colon E\to Q$ be a generically surjective morphism in a neighborhood
of the closure of the unit ball $D$, and let the currents $R,U$ and the Hefer forms
$H$ be defined as above. If $g$ for instance  is the smooth weight from Example~\ref{extva},
then for any holomorphic section $\phi$ of $Q$ ($r$-tuple of holomorphic functions)
we have the explicit  holomorphic  decomposition
\begin{equation}\label{decomp}
\phi(z) = f(z)\int H^1 U\phi\w g+ \int H^0 R\phi \w g, \qquad z\in D.
\end{equation}
\end{thm}

\begin{ex}
This  division formula is
non-trivial even when  $Z$ is empty so that $R=0$, and as an example let us compute it 
explicitly in the ball in $\C^2$.
Since $Z$ is empty, $U$ is smooth, and therefore we can use the weight $g$ from \eqref{r2}, and so the formula is
$\phi(z)=f(z)T\phi(z)$, where
\begin{multline*}
T\phi(z)=\frac{1}{(2\pi i)^2}\int_{|\zeta|=1} H^1_1 U_1\phi\w
\frac{\partial|\zeta|^2\w(\dbar\partial|\zeta|^2)^2}{(1-\bar\zeta\cdot z)^2}+\\
\frac{1}{2\pi i}\int_{|\zeta|=1} H^1_2 U_2\phi\w
\frac{\partial|\zeta|^2\w\dbar\partial|\zeta|^2}{1-\bar\zeta\cdot z}
+
\int_{|\zeta|<1} H^1_3U_3\phi.
\end{multline*}
If $f=f^1\otimes\epsilon_1+\ldots +f^r\otimes\epsilon_r$ and
$\sigma=\sigma_1\otimes\epsilon_1^*+\ldots+\sigma_r\otimes\epsilon_r^*$ as above,
and $\phi=\phi_1\epsilon_1+\ldots +\phi_r\epsilon_r$,  then
(suppressing the basis elements $\epsilon_1\w\ldots\w\epsilon_r$ and its dual
$\epsilon_1^*\w\ldots\w\epsilon_r^*$) we have that
$$
U_1\phi=\sigma\phi=\sum_1^r\phi_k\sigma_k,\quad
U_2\phi=\lsigma\w\dbar\sigma\phi=\sigma_1\w\ldots\w\sigma_r\w\sum_1^r\dbar\phi_k\sigma_k,
$$
and
$$
U_3\phi=\sum_{j=1}^r\dbar\sigma_j\otimes\epsilon_j\w\sigma_1\w\ldots\w\sigma_r\w\sum_1^r\dbar\phi_k\sigma_k.
$$
Next we have to compute $H^1_k$ for $k=1,2,3$. To begin with,  $H^1_1=I_{E_0}$ whereas 
 $H^1_2$ has to be a holomorphic solution to
$$
\delta_{\zeta-z} H^1_2=\delta_{f^1(\zeta)}\cdots\delta_{f^r(\zeta)}-\delta_{f^1(z)}\cdots\delta_{f^r(z)}.
$$
If $r=2$ one can take, e.g., 
$
H^1_2=\delta_{f^1}\delta_{h_2}-\delta_{h_1}\delta_{f^2(z)},
$
where $h=h_1\otimes\epsilon_1+\ldots +h_r\otimes \epsilon_r$ and $h_j$ are $(1,0)$-forms such that
$\delta_{\zeta-z} h_j=f^j(\zeta)-f^j(z)$.
Finally, $H^1_3$ has to solve (recall that $r$ is assumed to be odd;
in case $r$ is even, the first term on the right  should have a minus sign)
$$
\delta_{\zeta-z} H^1_3=H^1_2\delta_{f(\zeta)}-\delta_{f^1(z)}\cdots\delta_{f^r(z)}\delta_h,
$$
i.e.,
$H^1_3=(H^1_3)_1\otimes\delta_{\epsilon_1}+\ldots +(H^1_3)_r\otimes\delta_{\epsilon_r}$,
where
$$
\delta_{\zeta-z}(H^1_3)_k=H^1_2\delta_{f^k(\zeta)}-\delta_{f^1(z)}\cdots\delta_{f^r(z)}\delta_{h_k}.
$$
\end{ex}

\section{Various  applications}
In this section we illustrate the utility of  our new formula by presenting
some matrix variants of previously known results. In most cases, the
proofs are very similar the case $r=1$, so we only indicate them.

\subsection{A cohomological  duality result}
Let $f\colon E\to Q$ be generically surjective,  assume 
that $\codim Z=m-r+1$, and let $U$ and $R$ be the currents from the preceding section. 
We then know from \cite{A5} that 
$f\psi=\phi$ has a holomorphic  solution locally if and only if  $R\phi=R_{m-r+1}\phi=0$.  
Moreover, a solution $\psi$ is given by \eqref{soll}.

For degree reasons,   $R_{m-r+1}=\dbar U_{m-r+1}$, and we have a mapping
$$
G\phi\colon \xi\mapsto \int \dbar\xi \w u_{m-r+1}\phi
$$
for test forms $\xi$ such that $\dbar\xi=0$ in a neighborhood of $Z$.

\begin{prop} If $\codim Z=m-r+1$, then 
$G\phi=0$ if and only if $f\psi=\phi$ locally  has holomorphic solutions.
\end{prop}

In the case $r=1$ this  duality result  is proved in \cite{DS} and \cite{P1}.  
In the case that $m=n$, then $G\phi$  is the classical Grothendieck residue.
One can prove (see \cite{AW}) that $G$ as well as $R$ are  independent of the
Hermitian metrics on $E$ and $Q$, and essentially only depends on the sheaf $J=\Im (\O(E)\to\O(Q))$.

\begin{proof}
By Stokes' theorem it follows that $G\phi=0$ if $R_{m-r+1}\phi=0$, i.e.,
if $f\psi=\phi$ is locally solvable.
To prove the converse, we  mimick  the argument given in \cite{P1}
(the proof of Theorem 6.3.1). Clearly the statement is local, so let us fix a point on $Z$ that we may
assume is the origin.
After a suitable linear change of coordinates we may
assume that if  $W=\{|z'|, |z''|<1\}$, where  $z=(z',z'')\in\C^{n-(m-r+1)}\times\C^{m-r+1}$,
then $Z\cap W$ is contained in $\{|z''|<\delta\}$. Take $\chi=\chi'\chi''$, where
$\chi'=\chi'(z')$ has support in $|z'|<\delta$ and is identically $1$ for small $z'$, and
$\chi''$ is a cutoff function that is $1$ in a neighborhood $W$. Moreover, take
$$
s=\frac{\partial|\zeta'|^2}{2\pi i(|\zeta'|^2-\bar\zeta'\cdot z')} 
$$
in $W$ and extend it  outside $W$ so that $\delta_{\zeta-z}s\neq 0$ for 
$z$ close to the origin, and depends holomorphically of $z$ there.
Let $g=\chi -\dbar\chi\w s/\nabla_{\zeta-z}s$
as before. In $W$ then $g$ only depends on $z'$ so for degree reasons
$g_{\mu,\mu}$, $\mu=n-(m-r+1)$,   is $\dbar$-closed there. 
The obstruction for $T\phi$  being a solution to $f\psi=\phi$ is the residue term, cf., \eqref{decomp}
$$
\int H^0R_{m-r+1}\phi\w g_{\mu,\mu}.
$$
However, $R_{m-r+1}=\dbar U_{m-r+1}$ so an  integration by part gives 
$$
\int H^0 u_{m-r+1}\phi\w\dbar g_{\mu,\mu},
$$
which by assumption vanishes for  $z$ close to $0$.
This proves the statement.
\end{proof}

\subsection{A division problem  for smooth sections}

Let $f\colon E\to Q$ be a generically surjective morphism as before, and assume to begin
with that $\codim Z=m-r+1$. 
Let $\phi$ be a smooth section of $Q$ and assume that there is a smooth section of
$E$ such that $f\psi=\phi$. Arguing as in  the proof of Theorem~1.2 in \cite{A5}
it  follows that  $R\phi=R (f\psi)=(\dbar R')\psi=0$.
Let $\dbar^\alpha=\partial^\alpha/\partial \bar z^\alpha$.
Then $\dbar^\alpha\phi= f(\dbar^\alpha\psi)$, and it therefore follows that 
\begin{equation}\label{alag}
R(\dbar^\alpha \phi)=0
\end{equation}
for all $\alpha$.  For a general $f$, i.e., not necessarily such that $\codim Z=m-r+1$,
we  have the converse statement.

\begin{thm}
Suppose that $\phi$ smooth and assume that \eqref{alag} holds for all
$\alpha$.  Then $f\psi=\phi$ has a smooth solution.
\end{thm}

This was first proved for $r=1$ in \cite{A3}. We do not know any argument based on the
Koszul complex and successively solving of $\dbar$-equations as in the proof of
Proposition~\ref{locke}  above.
However, in \cite{JEB} is recently given a quite  simple proof based on a deep criterion for closedness
of ideals of smooth functions in terms of formal power series due to
Malgrange, \cite{Mal}.

If we replace $C^\infty$ with real-analytic functions $C^\omega$, then
the corresponding statement follows directly from the holomorphic
case, by embedding $X$ in the anti-diagonal
$\{(z,\bar z)\in\C^{2n};\ z\in X\}$.

\begin{cor}
Suppose that $\phi$ smooth and $\|\dbar^\alpha\phi\|\lesssim \det(ff^*)^{min(n,m-r+1)}$
for all $\alpha$. Then  $f\psi=\phi$ has a smooth solution.
\end{cor}

\begin{remark}
The corollary   can be seen as an extension of the Brian\c con-Skoda
theorem and follows by a standard estimate from the theorem.  

In the real-analytic case  it is easy to see that the size condition in the corollary is fulfilled
if $r=1$ and  $\phi$ belongs to the integral closure of the ideal $J=\E(f)$, i.e., if there
are functions $a_k \in J^k$ such that
$\phi^N +a_1\phi^{N-1}+\cdots +a_N=0.$
We do not know if the same is true in the smooth case.
\end{remark}

We also have an analogous result for lower regularity.

\begin{thm} Assume that $M$ is the order of the current $R$. There is a number $c_n$,
only depending on $n$, such that  if $\phi\in C^{c_n+2M+k}$ and
\eqref{alag} holds for all $|\alpha|\le c_n+M+k$,
then there is a section $\psi$ of $E$ of class $C^k$ such that
$f\psi=\phi$.
\end{thm}

Once we have the appropriate division formula these theorems follows in
the same way as for the case $r=1$ in \cite{A3}, and we therefore omit
the proofs.

\subsection{Matrix $H^p$-corona theorems}\label{mhp}
Suppose that $F(z)$ is a pointwise surjective  $r\times m$-matrix
of bounded holomorphic functions in a strictly pseudoconex domain $D$ and
assume furthermore that $|\det f(z)|\ge\delta>0$.
Then for each  $r$-tuple $\phi$ in $H^p(D)$, $p<\infty$, one can find an $m$-tuple
$\psi$ in $H^p$ such that $F\psi=\phi$.
This was proved in \cite{JH1} (and with a sharper estimate in \cite{A5},
see \cite{A5} also for a further discussion), by reducing it to
the case $r=1$ via the Fuhrmann trick, \cite{Fu}, and the case $r=1$
is known since long ago, see \cite{AC1} and the references  given there.
An explicit solution formula
in case $r=1$ is given in \cite{AC}, and copying the arguments  there, and
using the special choice of Hefer forms defined in Section~\ref{heferball}, (most likely) 
$$
T\phi= \int_{D}H^1 U\phi\w  g
$$
is such a solution in $H^p$ in the unit ball provided that
$$
g=\Big(\frac{1-\bar\zeta\cdot z}{1-|\zeta|^2}-\omega\Big)^{-\alpha},
$$
and $\alpha$ is large enough.

\subsection{Division formulas for $\dbar$-closed forms}

In \cite{A5} we  proved Brian\c con-Skoda type results also  for
$\dbar$-closed smooth $(p,q)$-forms, and even in this case we can provide
explicit representations of the solutions.
Again let $f\colon E\to Q$ be a generically surjective holomorphic morphism.
We want explicit expression for a $\dbar$-closed solution to
$f\psi=\phi$ provided that  $R\phi=0$.
Following \cite{A1} we now consider forms in $X_\zeta\times X_z$ with values
in the exterior algebra spanned by $T^*_{0,1}(X\times X)$ and
the $(1,0)$-forms $d\eta_1,\ldots,d\eta_n$, where  $\eta_k=\zeta_k-z_k$.
Then interior multiplication $\delta_\eta$ with
$\eta=\sum_1^n\eta_j(\partial/\partial \eta_j)$ has a meaning and
we can build up formulas pretty much as when $z$ is just considered as a
parameter. If we let 
$v=b/\nabla_\eta b$, where $b=(2\pi i)^{-1}\sum\partial|\eta|^2/|\eta|^2$,
and $\nabla_\eta=\delta_\eta-\dbar$, then
$
\nabla_\eta v=1-[\Delta],
$
where $[\Delta]$ is the $(n,n)$-current of integration over the diagonal
in $X\times X$. Let $g'=f(z)H^1U+H^0R$ as before  and let $g$ be the
form from Example~\ref{extva} but with all $d\zeta_k$ replaced by $d\eta_k$.
Then $\nabla_\eta (g'\w g)=0$ as before and therefore
at least formally we have 
$$
\nabla_\eta(v\w g'\w g)=g'\w g-[\Delta]I_Q,
$$
and for degree reasons thus
$$
\dbar (v\w g'\w g)_{n,n-1}= [\Delta]- (g'\w g)_{n,n}.
$$
Therefore,
$$
\phi(z)=\pm\int_\zeta v\w g\w g'\w\dbar\phi\pm\dbar_z\int_\zeta v\w g\w g'\w \phi
+\int g\w g'\w \phi, \quad z\in \U.
$$
Assuming that $\dbar\phi=0$, we get
\begin{multline}\label{brutt}
\phi(z)=\pm f(z)\dbar_z\int_\zeta H^1U(\zeta)\w v\w  g\w \phi+
f(z)\int_\zeta g\w H^1U\w\phi  +\\
\pm \dbar_z\int_\zeta v\w g\w  H^0R\w \phi +\int_\zeta g\w H^0R\w \phi.
\end{multline}

Since the integrals with $v$ are essentially  convolutions  of currents  and the
locally integrable functions $\zeta\mapsto \zeta_j/|\zeta|^{2k}$, $k\le n$,
they have  meaning.  One can prove \eqref{brutt} strictly by a suitable
approximation argument that we omit.

\smallskip
If $\phi$ is holomorphic, then we get back formula \eqref{decomp}.
Notice that $g$ and $g'$ are holomorphic in $z$ and therefore cannot contain
any component of positive degree in $d\bar z$. Therefore we have

\begin{prop}
If $\phi$ is a $\dbar$-closed $(p,q)$-form with values in $Q$, $q>0$,
such that $R\phi=0$, then  a $\dbar$-closed solution to
$f\psi=\phi$ is provided by the formula
$$
\psi(z)= \dbar_z\int_\zeta H^1U \w v\w  g\w \phi.
$$
\end{prop}

Notice that in general $\psi$ is not, and cannot be,  smooth.

\def\listing#1#2#3{{\sc #1}:\ {\it #2},\ #3.}

\end{document}